\documentclass[final,5p,times,twocolumn]{elsarticle}
\usepackage{amssymb,color,ulem}

\usepackage[usenames,dvipsnames]{xcolor}

\journal{European Journal of Control}
\usepackage{amsthm,verbatim}
\theoremstyle{definition}

\newtheorem{thm}{Theorem}[section]
\newtheorem{lem}{Lemma}[section]

\newcommand{\eps}{\varepsilon}
\begin{document}
\begin{frontmatter}
\title{A linear state feedback switching rule for global stabilization of switched nonlinear systems about a nonequilibrium point}
\author{Oleg Makarenkov}
\ead{makarenkov@utdallas.edu}
\address{Department of Mathematical Sciences, The University of Texas at Dallas 
800 West Campbell Road 
Richardson, TX 75080}
\begin{abstract} A switched equilibrium of a switched system of two subsystems is a such a point where the vector fields of the two subsystems point strictly towards one another. Using the concept of stable convex combination that was developed by Wicks-Peleties-DeCarlo (1998) for linear systems, Bolzern-Spinelli (2004)  offered a design of a state feedback switching rule that is capable to stabilize an affine switched system
to any switched equilibrium. 
The state feedback switching rule of Bolzern-Spinelli gives a nonlinear (quadratic) switching threshold passing through the switched equilibrium. In this paper we prove that the switching threshold (i.e. the associated switching rule) can be chosen linear, if each of the subsystems of the switched system under consideration are stable.

\end{abstract}
\begin{keyword} Switched system, switched equilibrium, global quadratic stabilization
\MSC[2010] 34H15  \sep 93D15  \sep 34A36  
\end{keyword}
\end{frontmatter}
\section{Introduction}\label{sec:int}
Using the concept of stable convex combination that was developed by Wicks et al  \cite{eur} for linear systems, Bolzern-Spinelli \cite{bol} offered a design of a state feedback switching rule that is capable to stabilize an affine switched system\footnote{Bolzern-Spinelli \cite{bol} actually considered a slightly more general case $\sigma:[0,\infty)\to\{1,...,m\}$, but in this paper we stick to just two discrete states.}
\begin{equation}\label{ls}
\dot x=A^\sigma x+b^\sigma, \quad x\in\mathbb{R}^n,\ \  \sigma\in\{-1,1\}
\end{equation}
to any point $x_0$ (called {\it switched equilibrium}) that satisfies 
\begin{equation}\label{se}
      \lambda \left(A^+ x_0+b^+\right)+(1-\lambda)\left(A^- x_0+b^-\right)=0,
\end{equation} 
for some $\lambda\in[0,1].$ If the matrix 
$\lambda A^++(1-\lambda)A^-$
is Hurwitz, 
then, according to Bolzern-Spinelli \cite{bol}, the switching signal $\sigma(x)$ can be defined as 
\begin{equation}\label{rule1} \def\arraystretch{1.2}
\begin{array}{l}   \sigma(x)={\rm arg}\min\limits_{i\in\{-1,1\}}\{V'(x)(A^i x+b^i)\}=\\
\qquad ={\rm sign}\left(V'(x)(A^- x+b^-)-V'(x)(A^+ x+b^+)\right),
\end{array}
\end{equation}
where
 $V$ is the quadratic Lyapunov function of the linear system
$$
\dot x= \lambda \left(A^+ x+b^+\right)+(1-\lambda)\left(A^- x+b^-\right).
$$
When $A^-=A^+$, the rule (\ref{rule1}) reduces to 
\begin{equation}\label{rule4}
\sigma(x)={\rm sign}\left(V'(x)b^--V'(x)b^+\right),
\end{equation}
whose switching threshold $\left\{x\in\mathbb{R}^n:V'(x)b^--V'(x)b^+\right\}\ni x_0$ is a hyperplane, but in general the state feedback switching rule (\ref{rule1}) gives a nonlinear switching threshold (quadratic surface) passing through the switched equilibrium $x_0.$

\vskip0.2cm

\noindent  In this paper we provide a wider class of switched systems (\ref{ls}) that can be stabilized to a switched equilibrium by a linear switching rule. Specifically, we
show that the nonlinear switching rule (\ref{rule1}) can always be replaced with the linear one
\begin{equation}\label{rulelin}
   \sigma(x)={\rm sign}\left<x-x_0,\left[V''(x_0)(A^-(x_0)+b^-)\right]^T\right>,
\end{equation}
 when
the subsystems $\dot x=A^+x$ and $\dot x=A^-x$  admit a common quadratic Lyapunov function. Here $V''(x_0)$ doesn't depend on $x_0$ because $V$ is assumed quadratic. We also note that (\ref{rulelin}) coincides with (\ref{rule4}) when $A^-=A^+.$

\vskip0.2cm

\noindent The paper is organized as follows. In the next section of the paper we discuss the main idea behind the switching rule (\ref{rule1}), which is based on construction of suitable sets $\Omega^-$ and $\Omega^+$, such that any switching  rule  $\sigma(x)$ with the property 
$$
\sigma(x)=\left\{\begin{array}{lll}
-1 & {\rm if} & x\in\Omega^-,\\
1 & {\rm if} & x\in\Omega^+,
\end{array}\right.
$$
stabilizes (\ref{ls}) to $x_0.$
In section~3 we prove our main result (Theorem~\ref{thm1}), which offers a linear state feedback switching rule to stabilize a nonlinear switched system
\begin{equation}\label{nsw}
\dot x = f^\sigma(x), \quad x\in\mathbb{R}^n,\ \ \sigma\in\{-1,1\},
\end{equation}
to a switched equilibrium $x_0$.
We recall that, according to Demidovich \cite[Ch.~IV, \S281]{dem}, nonlinear systems (\ref{nsw}) admit a common quadratic Lyapunov function, if  
the simmetrized derivative
$$
   f^\sigma_x(x)+\left[f^\sigma_x(x)\right]^T
$$
is uniformly negative definite uniformly in $x\in\mathbb{R}^n$, and $\sigma\in\{-1,1\},$ see also Pavlov et al \cite{pav}.
The switching rule (\ref{rule2}) proposed in Theorem~\ref{thm1} takes the form (\ref{rulelin}) when  switched system (\ref{nsw}) is affine. The main discovery used in Theorem~\ref{thm1} is that, for subsystems of (\ref{nsw}) that admit a common quadratic Lyapunov function, the boundaries of $\Omega^-$ and $\Omega^+$ are contained in ellipsoids that touch one another at the point $x_0,$ see Fig.~\ref{OmegaL}. The proof uses a standard Lyapunov stability theorem that is also implicitly used in Bolzern-Spinelli \cite{bol}. Specifically, we use a Lyapunov stability theorem for Filippov systems with smooth Lyapunov functions, which is a particular case of more general results available e.g. in 
Shevitz-Paden \cite{paden1} or M.-Aguilara-Garcia \cite{garcia}. But since deriving the required Lyapunov theorem (Theorem~\ref{th6}) from \cite{garcia,paden1} is not very straightforward (and since we didn't find the exact required theorem elsewhere in the literature), we  added a proof for completeness, that we placed in the Appendix section.

\vskip0.2cm

\noindent In section~4 we consider  an application of Theorem~\ref{thm1} to a model of boost converter and, for illustration purposes, also implement the Bolzern-Spinelli rule (\ref{rule1}) for the same model.  Some further discussion on when the switching rule (\ref{rulelin}) coincides with (\ref{rule1}) is carried out in the conclusions section. 

\section{The idea of Wicks et al \cite{eur} and Bolzern-Spinelli \cite{bol}}

\noindent Recall that $x_0$ is a switched equilibrium for the nonlinear switched system (\ref{nsw}), if there exists $\lambda_0\in[0,1]$ such that
\begin{equation}\label{barlambda}
      \lambda_0 f^-(x_0)+(1-\lambda_0)f^+(x_0)=0.
\end{equation} 

\noindent Assume that the equilibrium $x_0$ of the convex combination 
\begin{equation}\label{onesystem}
  \dot x= \lambda_0 f^-(x)+(1-\lambda_0)f^+(x).
\end{equation} is asymptotically stable and let $V$ be the respective Lyapunov function satisfying
\begin{equation}\label{nocommonV}
\begin{array}{l}
   V'(x)\left( \lambda_0 f^-(x)+(1-\lambda_0) f^+(x)\right)<0\quad\mbox{for all}\ x\not=x_0.
\end{array}
\end{equation}

\vskip0.2cm

\noindent The fundamental idea of Bolzern-Spinelli \cite{bol} (who extended Wicks et al \cite{eur} to affine linear systems) is that for (\ref{nsw}) to stabilize to $x_0$,  the  switching rule $\sigma(x)$  must take the value  $\sigma(x)=-1$  in the region
\begin{equation}\label{Omega-} 
   \Omega^-=\left\{x:V'(x)f^-(x)<0\right\}
\end{equation}
and the value $\sigma(x)=+1$ in the region
\begin{equation}\label{Omega+} 
   \Omega^+=\left\{x:V'(x)f^+(x)<0\right\}.
\end{equation}
\begin{figure}[t]\center
\includegraphics[scale=0.75]{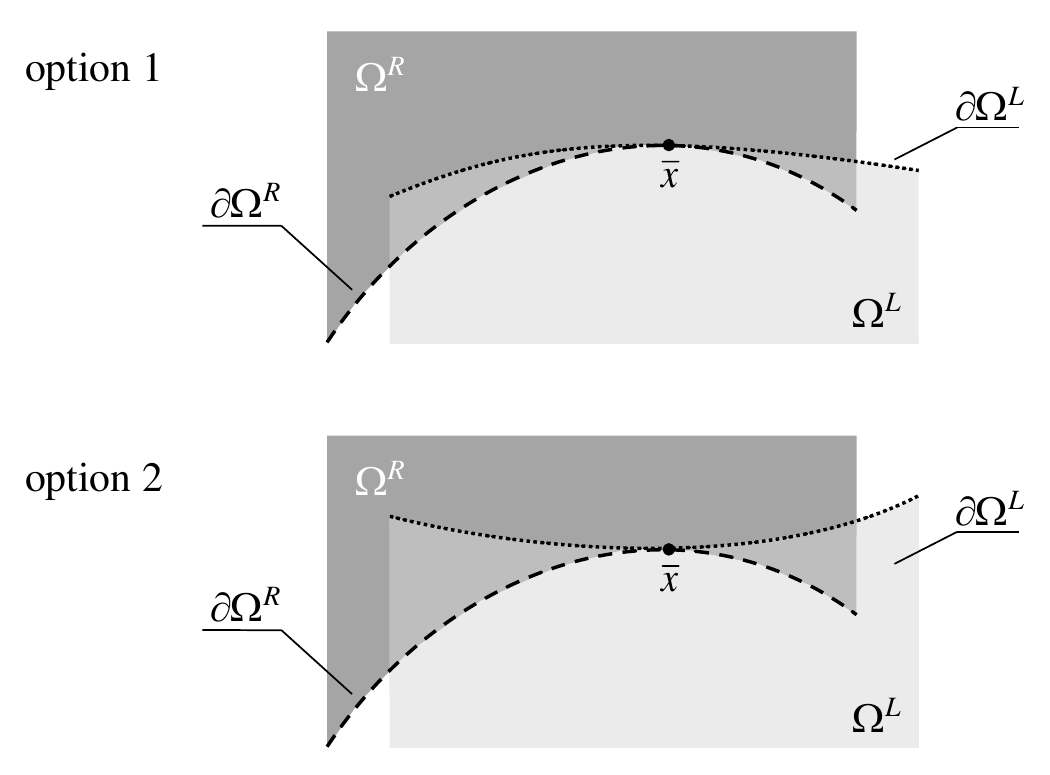}
\caption{\footnotesize Relative locations of sets $\Omega^L$ and $\Omega^R.$} \label{fig0}
\end{figure}
The following lemma discusses the geometry of the intersection $\Omega^-\cap\Omega^+$, in particular it clarifies that there are situations where one cannot draw a hyperlane in $\Omega^-\cap\Omega^+$ passing through $x_0$ (Fig.~\ref{fig0}a) and there are situations when one can (Fig.~\ref{fig0}b). The existence of a hyperplane in 
 $\Omega^-\cap\Omega^+$ passing through $x_0$ corresponds to the existence of a linear switching rule $\sigma(x)$ that stabilizes (\ref{nsw}) to $x_0.$
Therefore, what this paper will really prove in Section~3 is that it is Fig.~\ref{fig0}b which takes place when both of the subsystems of (\ref{nsw}) are stable.


\begin{lem}\label{ww} (ideas of  \cite{eur,bol}) Consider $f^-,\ f^+\in C^1(\mathbb{R}^n,\mathbb{R}^n)$.
Let $x_0$ be a switched equilibrium for the vector fields $f^-$ and $f^+$, i.e. (\ref{barlambda}) holds. Assume that the equilibrium $x_0$ of system (\ref{onesystem}) is asymptotically stable and the respective Lyapunov function $V\in C^1(\mathbb{R}^n,\mathbb{R})$ satisfies  (\ref{nocommonV}). Then, the sets $\Omega^-$ and $\Omega^+$ satisfy the properties:
\begin{itemize}
 \item [{\rm 1)}] $\Omega^-\cup\Omega^+\cup \{x_0\}=\mathbb{R}^n,$\ \ $\overline{\Omega^-}\cup\overline{\Omega^+}=\mathbb{R}^n,$ 
\item[{\rm 2)}] $\partial \Omega^-\backslash\{x_0\}\subset\Omega^+$,\ $\partial \Omega^+\backslash\{x_0\}\subset\Omega^-$,
\item[{\rm 3)}] $x_0\in\partial \Omega^-,$  $x_0\in\partial \Omega^+.$
\end{itemize}
\end{lem}


\noindent{\bf Proof.} {\bf Part 1.} Follows directly from (\ref{nocommonV}). 
\vskip0.2cm

\noindent {\bf Part 2.} Consider $x\in\partial \Omega^-$. Then $x\not\in\Omega^-$ because $\Omega^-$ is open. Then $x\in\overline{\Omega^+}$ by Part~1.  The property $\partial \Omega^+\subset \overline{\Omega^-}$ can be proved by analogy.

\vskip0.2cm

\noindent {\bf Part 3.} It is sufficient to show that $V'(x_0)=0$. To observe this, fix an arbitrary $j\in\overline{1,n}$ and consider the vector $\xi^j\in\mathbb{R}^n$ defined as $\xi_i^j=0,$ $i\not= j$, and $\xi_j^j=1.$ Since $V(x)>0$, $x\not=x_0$, we have
$$
\begin{array}{l}
   0<V(x_0+k \xi^j)-V(x_0)=V'(x_0+k_*\xi^j)\xi^j\cdot k=\\
   \qquad =\frac{\partial V}{\partial x_j}(x_0+k_*\xi^j)k,\\
0<V(x_0-k \xi^j)-V(x_0)=-V'(x_0-k_{**}\xi^j)\xi^j\cdot k=\\
\qquad=-\frac{\partial V}{\partial x_j}(x_0-k_{**}\xi^j)k,
\end{array}
$$
for any $k>0$ and for some $k_*,k_{**}\in[0,k]$ (that depend on $k$). Passing to the limit as $k\to 0$, one gets $\frac{\partial V}{\partial x_j}(x_0)=0$. 
\vskip0.2cm

\noindent The proof of the lemma is complete.\qed

\section{The main result}

\noindent In this section we assume that the switched equilibrium $x_0$ admits a common quadratic Lyapunov function 
$$
V(x)=(x-x_0)^TP(x-x_0)
$$ with respect to each of the two systems
\begin{equation}\label{asystems}
   \dot x=f^-(x)-f^-(x_0)\quad\mbox{and}\quad  \dot x=f^+(x)-f^+(x_0),
\end{equation}
where $P$ is an $n\times n$ symmetric matrix and the following standard properties hold:
\begin{equation}\label{acommonV}
\begin{array}{l}
   V'(x)\left(f^-(x)-f^-(x_0)\right)\le -\alpha\|x-x_0\|^2,\\    V'(x)\left(f^+(x)-f^+(x_0)\right)\le -\alpha\|x-x_0\|^2,
\end{array}
\end{equation}
for some fixed constant $\alpha>0$.

\vskip0.2cm

\begin{thm} \label{thm1} \it
Consider $f^-,\ f^+\in C^1(\mathbb{R}^n,\mathbb{R}^n)$.
Let $x_0$ be a switched equilibrium for the vector fields $f^+$ and $f^-$, i.e. (\ref{barlambda}) holds. 
Assume that the systems  of (\ref{asystems}) admit a common quadratic Lyapunov function $V\in C^2(\mathbb{R}^n,\mathbb{R})$ that satisfies (\ref{acommonV}). Then the switching signal
\begin{equation}\label{rule2}
  \sigma(x)={\rm sign}\left<x-x_0,\left[V''(x_0)f^-(x_0)\right]^T\right>
\end{equation}
makes $x_0$ quadratically globally stable switched equilibrium of switched system (\ref{nsw}).
\end{thm}

\noindent Note that rule (\ref{rule2}) takes the form (\ref{rulelin}) when the nonlinear switched system (\ref{nsw}) takes the form (\ref{ls}). Also, using (\ref{barlambda}) the switching rule (\ref{rule2}) can be rewritten as
$$
  \sigma(x)={\rm sign}\left<x-x_0,\left[V''(x_0)\left(f^-(x_0)-f^+(x_0)\right)\right]^T\right>.
$$

\vskip0.2cm

\noindent In order to prove the theorem, we 
introduce two sets
$$
\begin{array}{l}
    \Omega^-_\alpha=\left\{x\in\mathbb{R}^n:-\alpha\|x-x_0\|^2+V'(x)f^-(x_0)<0\right\},\\
    \Omega^+_\alpha=\left\{x\in\mathbb{R}^n:-\alpha\|x- x_0\|^2+V'(x)f^+(x_0)<0\right\}
\end{array}
$$
and establish the following lemma about the relative properties of the sets $\Omega^i_\alpha$ and $\Omega^i$ as introduced in (\ref{Omega-})-(\ref{Omega+}).

\begin{lem}\label{wwa} Assume that the conditions of Theorem~\ref{thm1} hold. Then $\Omega^-_\alpha$ and $\Omega^+_\alpha$ verify the following properties:
\begin{itemize}
\item[1)] $\Omega^-\supset \Omega^-_\alpha,$\ \  $\Omega^+\supset \Omega^+_\alpha,$
\item[2)] $x_0\in\partial \Omega^-_\alpha,$ \ \ $x_0\in\partial\Omega^+_\alpha$,
\item[3)] both $\partial\Omega^-_\alpha$ and $\partial\Omega^+_\alpha$ are ellipsoids,
\item[4)] hyperplane $\sigma(x)=0$ is tangent to both $\Omega^-_\alpha$ and $\Omega^+_\alpha$ at $x_0,$
\item[5)] $\Omega^-_\alpha\subset \left\{x:\sigma(x)<0\right\},$ \ $\Omega^+_\alpha\subset \left\{x:\sigma(x)>0\right\}.$
\end{itemize}
\end{lem}

\begin{figure}[t]\center
\includegraphics[scale=0.75]{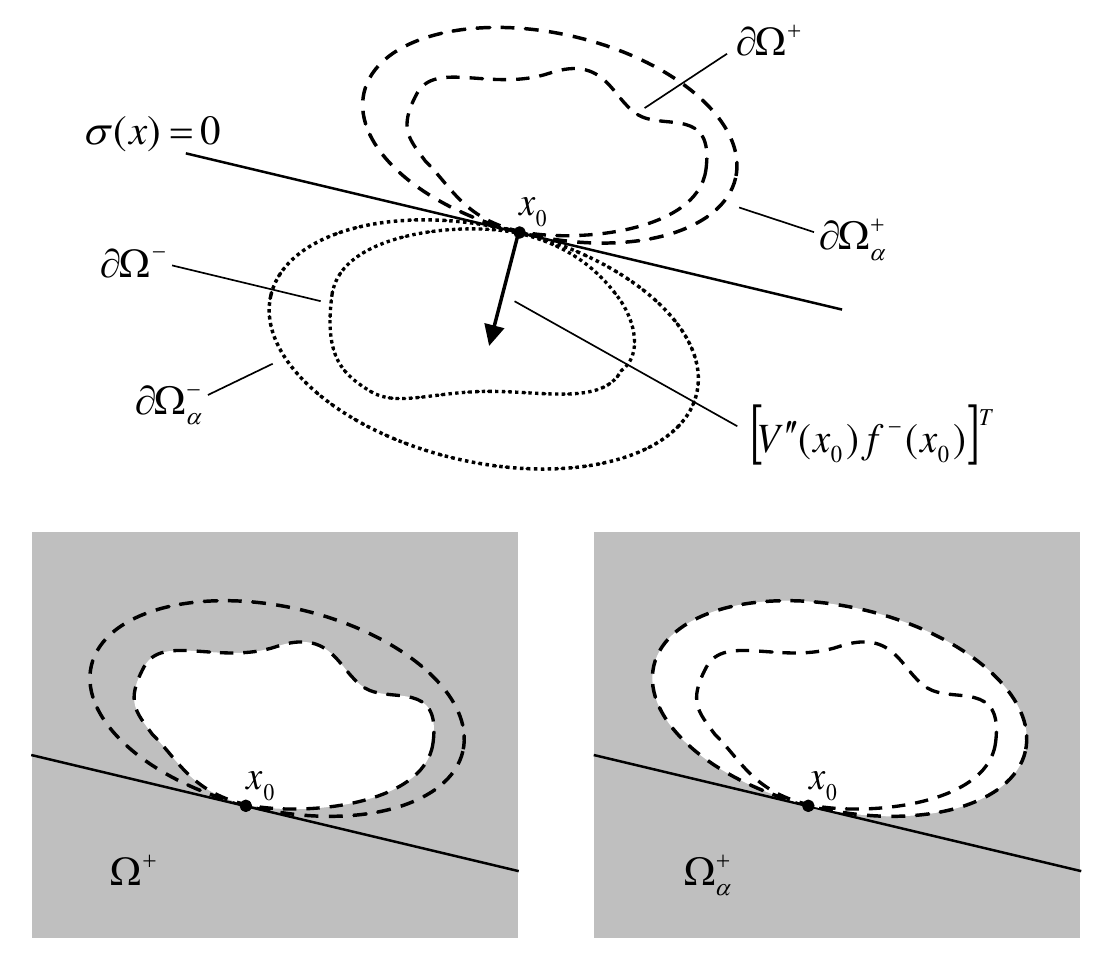}
\caption{\footnotesize Top figure: Locations of the boundaries of $\Omega^+,$ $\Omega^+_\alpha$, $\Omega^-,$ $\Omega^-_\alpha$ with respect to the hyperplane $\sigma(x)=0$ and with respect to each other. 
Bottom figures: The sets $\Omega^+$ and $\Omega^+_\alpha$ (grey regions).} \label{OmegaL}
\end{figure}

\noindent The notations and statements of Lemma~\ref{wwa} are illustrated at Fig.~\ref{OmegaL}.

\vskip0.2cm

\noindent {\bf Proof.} {\bf Part 1.} Let $x\in\Omega_\alpha^-$. Then
$$
\begin{array}{l}
  V'(x)f^-(x)=V'(x)(f^-(x)-f^-(x_0))+V'(x)f^-(x_0)\le \\
   \le -\alpha\|x-x_0\|^2+V'(x)f^-(x_0)<0.
\end{array}
$$
Therefore, $x\in\Omega^-$. The proof for  $\Omega_\alpha^+$ and $\Omega_\alpha^+$ is analogous.

\vskip0.2cm

\noindent {\bf Part 2.} Follows from $V'(x_0)=0$ established in the proof of Part~3 of Lemma~\ref{ww}.

\vskip0.2cm

\noindent {\bf Part 3.} We execute the proof for $x_0=0$. The proof in the general case doesn't differ. The change of the coordinates $y=x-\Delta$ transforms the equation 
$$
  -\alpha\|x-x_0\|^2+V'(x)f^-(x_0)=0
$$
into 
$$
   -\alpha\|y\|^2-2\alpha\left<\Delta,y\right>+2\left<Pf^-(0),y\right>-\alpha\|\Delta\|^2+2\left<\Delta,Pf^-(0)\right>=0.
$$
If $\Delta=\frac{Pf^-(0)}{\alpha},$ then we further get
$$
   -\alpha\|y\|^2-\frac{1}{\alpha}\|Pf^-(0)\|^2+\frac{2}{\alpha}\|Pf^-(0)\|^2=0,
$$
which is the equation of ellipsoid centered at 0 and radius $\frac{1}{\alpha^2}\|Pf^-(0)\|^2.$
\vskip0.2cm

\noindent The proof for $\partial\Omega_\alpha^+$ is analogous.

\vskip0.2cm

\noindent {\bf Part 4.} This follows from the equality
$$
   \left.\frac{d}{dx}\left(-\alpha\|x-x_0\|^2+V'(x)f^-(x_0)\right)\right|_{x=x_0}=V''(x_0)f^-(x_0).
$$
and the property (\ref {barlambda}) of switched equilibrium.

\vskip0.2cm

\noindent {\bf Part 5.} Let $H(x)=-\alpha\|x-x_0\|^2+V'(x)f^-(x_0).$ The interior of the ellipsoid $\partial \Omega_\alpha^-$ corresponds to $H(x)>0$. Therefore, the exterior of the ellipsoid $\partial \Omega_\alpha^-$ (which, by definition, coincides with the set $\Omega_\alpha^-$) corresponds to $H(x)<0.$ This proves  the statement of Part~5 for $\Omega_\alpha^-.$
Since $(1-\lambda_0)f^+(x_0)=-\lambda f^-(x_0)$ by (\ref{barlambda}), the proof for $\Omega_\alpha^+$ follows same lines.

\vskip0.2cm

\noindent The proof of the lemma is complete. \qed



\vskip0.2cm

\noindent The proof of our main result uses the following  Lyapunov stability theorem for discontinuous systems with smooth Lyapunov functions, which is implicitly used in \cite{eur,bol}.

\begin{thm}\label{th6}{\bf (Lyapunov stability theorem for discontinuous systems with smooth Lyapunov functions)} {\rm (similar to \cite[Theorem 3.1]{paden1}, \cite[Theorem 2.3]{garcia})} Consider a system of differential equations with discontinuous right-hand-side
\begin{equation}\label{fil}
 \hskip-0.65cm  \dot x= g(x),\ \ {\rm with} \ \ g(x)=\left\{\begin{array}{l}
g^+(x),\ {\rm if}\ H(x)>0,\\
g^-(x),\ {\rm if}\ H(x)<0,
\end{array}\right.\ \ x\in\mathbb{R}^n, 
\end{equation}
where $g^-,$ $g^+$, and $H$ are $C^1$-functions. Consider $x_0\in\mathbb{R}^n$ satisfying $H(x_0)=0.$ Let $V$ be a $C^1$-smooth Lyapunov function with $V(x_0)=0$ and $V(x)\not=0$ for $x\not=x_0.$ 
 Consider a piecewise continuous strictly positive for $x\not=x_0$ scalar function $x\mapsto w(x)$   such that for any $\rho>0$ there exists $\eps>0$  for which $w(x)\ge \eps$ as long as $\|x-x_0\|\ge\rho.$ If 
$$
    V'(x)\xi\le-w(x)\quad\mbox{for any}\ \xi\in K[g](x),\ \mbox{and any}\ x\not=x_0,
$$
then $x_0$ is an asymptotically globally stable stationary point of (\ref{fil}). Here $K[g](x)$ stays for convexification of the discontinuous function $g$ at $x$, see e.g. Shevitz-Paden \cite{paden1}.
\end{thm}

\noindent The proof of Theorem~\ref{th6} is given in Appendix.

\vskip0.2cm

\noindent {\bf Proof of Theorem~\ref{thm1}.} We will show that the conditions of Theorem~\ref{th6} hold with 
$$w(x)=\left\{\begin{array}{ll}
-V'(x)f^-(x), &  \sigma(x)<0,\\
-\max\{V'(x)f^-(x),V'(x)f^+(x)\}, & \sigma(x)=0, \\
-V'(x)f^+(x), & \sigma(x)>0.
\end{array}\right.
$$
If $x\in \overline{D^-}\backslash\{x_0\}$, then $   x\in \Omega^-_\alpha\subset \Omega^-$  by statements 5 and 1 of Lemma~\ref{wwa}, which implies $w(x)>0$. Analogously, $w(x)>0$, if $x\in \overline{D^+}\backslash\{x_0\}.$ This implies that $\max\limits_{x:\|x-x_0\|=\rho}w(x)$ is a positive function of $\rho$ that approaches 0 as $\rho\to 0.$ 

\vskip0.2cm

\noindent Since 
$K[f](x)=\{f^-(x)\},$ when $\sigma(x)<0$, and $K[f](x)=\{f^+(x)\},$ when $\sigma(x)>0$, then condition  $V'(x)\xi\le -w(x)$ of Theorem~\ref{th6} holds for $\sigma(x)\not=0.$

\vskip0.3cm

\noindent Consider $\sigma(x)=0.$ Then each $\xi\in K[f](x)$ has the form
$
   \xi=\lambda f^-(x)+(1-\lambda)f^+(x),
$ where $\lambda$ is a constant from the interval $[0,1].$ We have
\begin{eqnarray*}
   V'(x)\xi&=&\lambda V'(x)f^-(x)+(1-\lambda)V'(x)f^+(x)\le\\
   && \le\max\{V'(x)f^-(x),V'(x)f^+(x)\}=-w(x),
\end{eqnarray*}
that completes the proof of the theorem. \qed

\section{Application to a model of boost converter}

\begin{figure}[t]\center
\includegraphics[scale=0.75]{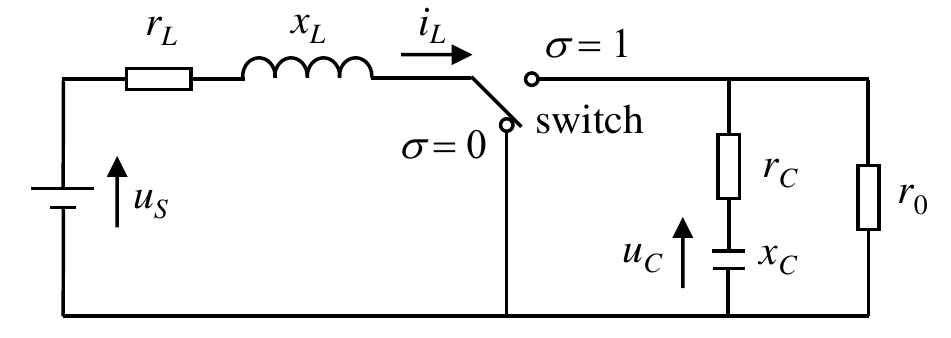}
\caption{\footnotesize Boost converter from Fribourg-Soulat \cite{boost2} and Beccuti et al \cite{boost1}.} \label{buck}
\end{figure}

\noindent Consider a dc-dc boost converter of Fig.~\ref{buck} with a switching feedback $\sigma(x).$ Denoting 
the inductor current $i_L$ by $x_1$ and the capacitor voltage $u_C$ by $x_2$, the differential equations of the converter read as (see e.g. Fribourg-Soulat \cite{boost2}, Beccuti et al \cite{boost1}) 
\begin{equation}\label{ex}\def\arraystretch{1.5}
 \hskip-0.6cm  \dot x=\left(\begin{array}{cc} -\frac{r_L}{x_L} & -\frac{r_0}{x_L(r_0+r_C)}\sigma \\
     \frac{r_0}{x_C(r_0+r_C)}\sigma & -\frac{1}{x_C(r_0+r_C)}\end{array}\right)x+\left(\begin{array}{c} \frac{u_s}{x_L} \\ 0\end{array}\right),\ \ \sigma\in\{0,1\},
\end{equation}
Let us view the right-hand-side of (\ref{ex}) with $\sigma=0$ and $\sigma=1$ as $f^-(x)$ and $f^+(x)$ respectively.
The equation  (\ref{barlambda}) for switched equilibrium $x_0$ yields
\begin{equation}\label{ex1}\def\arraystretch{1.5}
\hskip-0.6cm\begin{array}{l}
    -r_Lx_{01}+u_s-(1-\lambda_0)\frac{r_0 r_C}{r_0+r_C}x_{01}-(1-\lambda_0)\frac{r_0}{r_0+r_C}x_{02}=0,\\
    -x_{02}+(1-\lambda_0)r_0x_{01}=0,
\end{array}
\end{equation}
which can be solved for $(x_{01},\lambda_0)$ when the reference voltage $x_{02}$ is fixed.
 The conditions of Theorem~\ref{thm1} hold with the Lyapunov function
$$
   V(x)=\frac{1}{2x_C}(x_1-x_{01})^2+\frac{1}{2x_L}(x_2-x_{02})^2.
$$
 Therefore, 
$$
V''(x_0)f^-(x_0)=\left(\frac{1}{x_C}\left(-\frac{r_L}{x_L}x_{01}+\frac{v_s}{x_L}\right),\frac{1}{x_L}\left(-\frac{1}{x_C(r_0+r_C)}x_{02}\right)\right),
$$
which transpose will be denoted by $n.$
Plugging $n$ into (\ref{rule2}), we conclude that any point $x_0$ that satisfies the switched equilibrium condition (\ref{ex1}) with $\lambda_0\in(0,1)$, can be stabilized using the switching rule
\begin{equation}\label{rule3}
  \sigma(x)=\left\{\begin{array}{lll} 1, &  {\rm if} & (x-x_d)n>0,\\ 
  0, &  {\rm if} & (x-x_d)n<0.
  \end{array}\right.
\end{equation}

\begin{figure}[t]\center
\vskip-0.15cm
\includegraphics[scale=0.6]{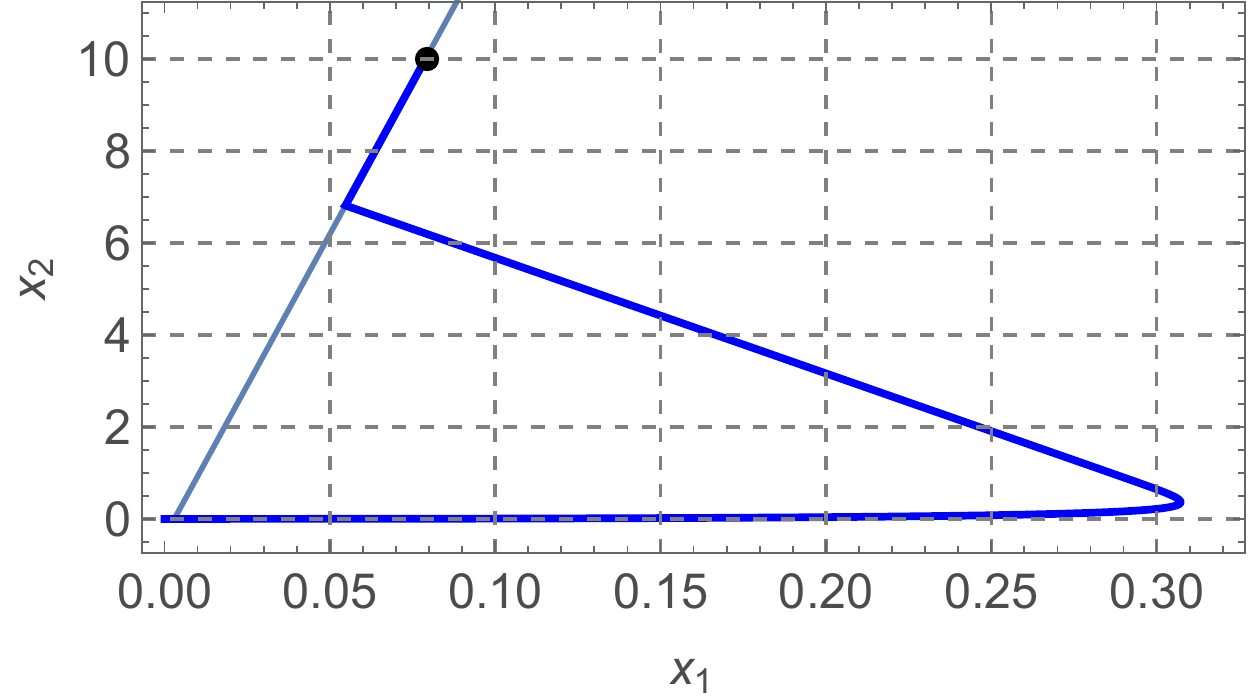}\\ 
\vskip0.4cm
\includegraphics[scale=0.6]{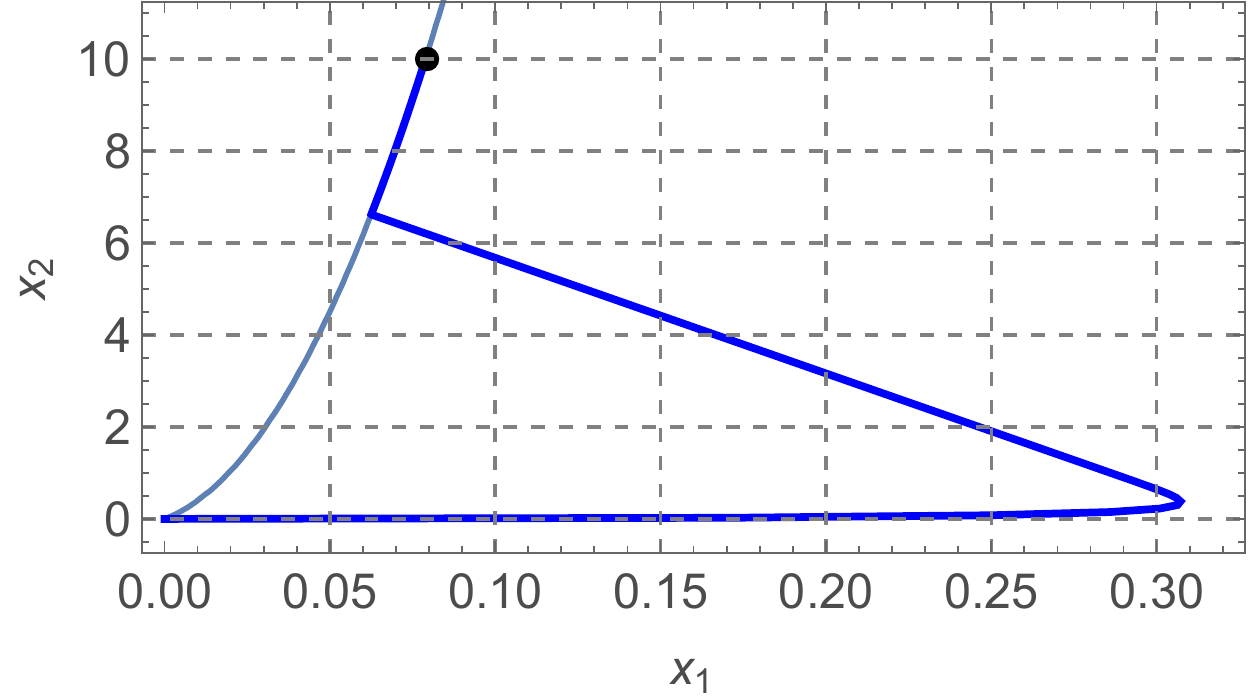}
\caption{\footnotesize The solution (bold curve) of switched system (\ref{ex1}) with the initial condition $x(0)=0$, the parameters (\ref{para}), and the switching signal $\sigma(x)$ given by (\ref{rule3}) (top figure) and by (\ref{spi}) (bottom figure). The thin curve is the switching manifold $\sigma(x)=0$ and the bold point is the switched equilibrium $x_0$. } \label{simfig}
\end{figure}

\noindent An implementation of switching rule (\ref{rule3}) with the parameters 
\begin{equation}\label{para}
  \hskip-0.7cm  r_L=20,\ r_C=5,\ x_L=600,\ x_C=70,\ r_0=200,\ u_s=8,
\end{equation}
and the reference voltage $x_{02}=10$ (which, when plugged into (\ref{ex1}), yields $x_{01}=0.79$ and $\lambda_0=0.367$ as one of the two possible solutions) 
is given in Fig.~\ref{simfig} (top). 

\vskip0.2cm

\noindent For comparison, Fig.~\ref{simfig} (bottom) shows stabilization of (\ref{ex1}) to the switched equilibrium $x_0=(0.79,10)$ using the switching rule (\ref{rule1}), that can be shown to simplify to 
\begin{equation}\label{spi}
  \sigma(x)={\rm sign}\left(r_Lr_Cx_1^2-(x_{01}r_Lr_C-x_{02})x_1-x_{01}x_2\right).
\end{equation}

\noindent The parameters (\ref{para}) are slightly artificial, but similar to Fig.~\ref{simfig} simulations are achieved in the case of more realistic parameters 
e.g. taken from \cite{boost1}, \cite{boost2}, or \cite{boost3}. The parameters (\ref{para}) are chosen in such a way that the nonlinear behavior of the Bolzern-Spinelli rule (\ref{spi}) is clearly seen in Fig.~\ref{simfig} (bottom). The top and bottom figures of Fig.~\ref{simfig} turn out to be indistinguishable (on the screen) for the parameters from \cite{boost1,boost2,boost3}.

\section{Conclusions} \noindent In this paper we showed that the switching rule (\ref{rule1}) of  Bolzern-Spinelli \cite{bol} for quadratic stabilization of a switched equilibrium $x_0$ of switched system (\ref{ls})  can be replaced by a linear switching rule (\ref{rulelin}) when the subsystems of (\ref{ls}) admit a common quadratic Lyapunov function. Moreover, our main result (Theorem~\ref{thm1}) applies to nonlinear switched systems  (\ref{nsw}) complimenting  the work by Mastellone et al \cite{spong} that proposes a nonlinear extension of  Bolzern-Spinelli \cite{bol} in the case where the subsystems of (\ref{nsw}) are shifts of one another (at the same time, the work \cite{spong} addresses the case of an arbitrary number of  subsystems, while the present paper focuses on just two subsystems).

\vskip0.2cm

\noindent We would like to note that seemingly nonlinear switching rule (\ref{rule1}) of Bolzern-Spinelli \cite{bol} simplifies to linear in wide classes of particular applications, e.g. in applications to buck converters (see e.g. Lu et al \cite{lu}), where $A^+=A^-$ in (\ref{ls}), or in applications to boost converters of Fig.~\ref{buck} with neglected resistance $r_C$ of the capacitor (see e.g. Schild et al \cite{sch}).  Still, the switching rule (\ref{rule1}) stays nonlinear in some other classes of applications, e.g. in more general boost converters such as the one of Fig.~\ref{buck} or its further extensions (see Gupta-Patra \cite{boost3} and references therein). In these classes of applications the linear switching rules (\ref{rulelin}) and (\ref{rule2}) proposed in this paper may simplify the engineering implementation of the feedback control.

\section{Appendix: Lyapunov stability theorem for  discontinuous systems with smooth Lyapunov functions}

\noindent {\bf Proof of Theorem~\ref{th6}.}  Let $x$ be a Filippov solution of (\ref{fil}), see e.g. Shevitz-Paden \cite{paden1}.  We pick  $\rho>0$ and prove that $x(t)\in \overline{W_\rho}$ beginning  some $t=t_\rho,$ where
$
{W_\rho}=\{x\in\mathbb{R}^n:V(x)< \rho\}.
$ 

\vskip0.2cm

\noindent {\bf Step 1.} Let $r>0$ be such a constant that $x(0)\in \partial W_r.$ We claim that $x(t)\in W_r$ for all $t>0.$
We prove by contradiction, i.e. assume that $x(\tau)\not\in W_r$ for some $\tau>0.$ Without loss of generality we can assume that $x([0,\tau])\subset W,$ where $W$ is an open neighborhood of $\overline{W_r}$, such that $w(x)$ is strictly positive in $W\backslash\{x_0\}.$ For the function
$
  v(t)=V(x(t))
$
we have 
\begin{equation}\label{ftfftf}
    v(0)=r\quad\mbox{and}\quad v(\tau)\ge r.
\end{equation}

\noindent {\bf Step 1.1} We claim that $v(t)> r/2$ for all $t\in[0,\tau]$. Indeed, if the latter is wrong, then defining
$
   s=\max\left\{t\in[0,\tau]:v(t)\le r/2\right\},
   $
one gets
\begin{equation}\label{tc}
  \hskip-0.7cm  v(s)=r/2,\ v(\tau)=r, \ v(t)\in \left[r/2,r\right],\ \mbox{for any}\ t\in[s,\tau].
\end{equation}
In particular, $x(t)\not=x_0$ for all $t\in[s,\tau]$ and, therefore,
$$
   v'(t)=V'(x(t))\xi<0,$$ for some $\xi\in K[f](x(t))$ and almost any  $t\in[s,\tau].$
This contradicts (\ref{tc}) and proves that $v(t)>r/2$ for all $t\in[0,\tau].$

\vskip0.2cm

\noindent {\bf Step 1.2}  Step 1.1 implies that 
$
   x(t)\not=0,$  for any $t\in[0,\tau],$  and, as a  consequence, $$v'(t)<0,\  \mbox{for any}\ t\in[0,\tau],
$$
which contradicts (\ref{ftfftf}) and completes the proof of the fact that $x(t)\in W_r$ for all $t>0.$

\vskip0.2cm

\noindent {\bf Step 2.} Let us show that $x(t)$ reaches $\overline{W_\rho}$ at some time moment. Assume that $x(t)$ never reaches $\overline{W_\rho}$. Then 
$$
   v'(t)=V'(x(t))\xi<-w(x(t)),
   $$
for some  $\xi\in K[f](x(t))$ and almost any $t>0.$
The definition of function $w$ implies that 
$
  w_{\min}=\min\{w(x), \ x\in \overline{W_r}\backslash W_\rho\}>0. 
$
Therefore,
$$
   v(t)=v(0)+\int_0^t v'(t)dt<v(0)-w_{\min}t
$$
and $v(t)$ becomes negative, if $x(t)$ never reaches $\overline{W_\rho}$. Since $\rho\in(0,r)$ was chosen arbitrary, our conclusion implies that $x(t)\to x_0$ as $t\to\infty.$ 

\vskip0.2cm

\noindent The proof of the theorem is complete.
\qed

\section*{Acknowledgements}
The research was supported by NSF Grant CMMI-1436856.

\bibliographystyle{plain}

\end{document}